\documentclass[10pt]{article}
\usepackage{}
\usepackage{indentfirst}
\usepackage{amsfonts}
\usepackage{amsfonts}
\usepackage{amsmath}
\allowdisplaybreaks[4]
\usepackage{bbm}
\usepackage{subfigure}
\usepackage{mathrsfs}
\usepackage{graphicx,amsmath,amsfonts,amssymb,amscd,amsthm}
\usepackage[ urlcolor=blue, colorlinks,linkcolor=black,anchorcolor=blue,citecolor=blue]{hyperref}
\newtheoremstyle{kai}
{10pt} {10pt} {} {} {\bfseries} {.} {.5em} {}
\def\EquationsBySection{\def\theequation
{\thesection.\arabic{equation}}%
\@addtoreset{equation}{section}}

\newcommand\old[1]{}
\newcommand{\pend}{\hfill \thicklines \framebox(6.6,6.6)[l]{}}
\renewenvironment{proof}{\noindent {\it  Proof.} \rm}{\pend}
\newtheorem{theorem}{Theorem}[section]
\newtheorem{lemma}{Lemma}[section]

\newtheorem{remark}{Remark}[section]
\newtheorem{definition}{Definition}[section]

\newtheorem{assumption}{Assumption}[section]

 \EquationsBySection \makeatother

\begin{document}

\pagestyle{plain}
\title
{\bf Almost sure exponential stability of multi-dimensional stochastic systems driven by a symmetric $\alpha$-stable process}

\author{\small Xuekang Zhang$^{1}$\footnote{Corresponding author. E-mail: xkzhang@ahpu.edu.cn
(X. Zhang),\,\, gjshen@163.com(G. Shen).
}, Hongjian Liu$^{1}$, Guangjun Shen$^{2}
$
 \\
\footnotesize{$^1$ School of Mathematics-Physics and Finance, Anhui Polytechnic University, Wuhu,  241000, China;} \\
\footnotesize{$^2$ Department of Mathematics, Anhui Normal University, Wuhu, 241002, China.}
}

\date{}
\maketitle \indent {\bf Abstract.} This article investigates the stability of a class of multidimensional stochastic systems driven by a symmetric $\alpha$-stable process with $1<\alpha<2$. These driving processes have infinite second moments, which prevents the direct application of standard second-moment methods and may lead to divergent jump integrals in the generators of quadratic Lyapunov functions. We establish criteria for $p$th-moment exponential stability and almost sure exponential stability. A numerical example illustrates the theoretical results.


\section{Introduction}

Dynamical systems arising in engineering, biology, and finance are inevitably affected by environmental noise, which is frequently modeled by Brownian motion. Stochastic differential equations (SDEs) driven by Brownian motion have been widely used to model dynamical systems in various fields. Stability analysis is a fundamental topic in the study of stochastic systems, as it provides insight into their long-term behavior under random perturbations. The stability of stochastic systems driven by Brownian motion has been extensively studied (see, e.g., Appleby et al. \cite{Appleby2008}, Dong and Mao \cite{Dong2022}, Has'minskii \cite{Hasminskii1981}, Hu et al. \cite{Hu2023}, Mao \cite{Mao1994}, Xu et al. \cite{Xu2026}, Yin et al. \cite{Yin2011}, Yu et al. \cite{Yu2026}).  Among various notions of stability, exponential stability is particularly important because it quantifies the exponential rate of decay of system solutions. In particular, exponential stability has also been extensively studied (see, e.g., Mao \cite{mao2016}, Tang et al. \cite{Tang2026}, Yang et al. \cite{Yang2024}, Zhang and Zhu \cite{Zhang2024}, Zhao and Deng \cite{Zhao2024}).    For instance, Mao \cite{Mao2007} studied $p$th moment and almost sure exponential stability for multi-dimensional nonlinear SDEs driven by Brownian motion. Xu et al. \cite{Xu2024}  employed improved Razumikhin-type theorems to analyze the globally weakly stochastically exponential stability of stochastic nonlinear delay systems subject to multiple periodic impulses. Zhang et al. \cite{ZhangHu2025}  proved $p$th moment exponential stability of stochastic highly nonlinear delayed systems with regime-switching diffusion under asynchronous intermittent control via an auxiliary-timer Lyapunov functional and Dupire's functional It\^o's formula.

However, many real-world phenomena exhibit heavy-tailed behavior (see, e.g., Applebaum \cite{Applebaum2009}, Barabasi \cite{Barabasi2005}, Sasaki et al. \cite{Sasaki2017} ).
 Such phenomena motivate the use of non-Gaussian models with heavy-tailed distributions. The symmetric $\alpha$-stable (S$\alpha$S) distribution has been widely used
to model a variety of heavy-tailed phenomena with applications
including synthetic aperture radar (Banerjee et al. \cite{Banerjee1999}), electromagnetic
interference (Ilow and Hatzinakos \cite{Ilow1998}), wireless sensor and ad hoc networks (Gulati et al. \cite{Gulati2010}), and shallow water acoustic channel (Zha and Qiu \cite{Zha2006}). S$\alpha$S distributions contain the Gaussian distribution as the case $\alpha=2$. Moreover, the generalized central limit theorem states that if a non-degenerate limit exists for suitably normalized sums of independent identically distributed random variables, then this limit must be an $\alpha$-stable distribution. For $\alpha\in(1,2)$, an S$\alpha$S random variable has a finite absolute moment of order $p$ if and only if $p\in(0,\alpha)$, whereas its variance is infinite. An S$\alpha$S process has stationary independent increments and S$\alpha$S marginal distributions. For $\alpha\in(1,2)$, an S$\alpha$S process has jump-type sample paths and heavy-tailed marginal distributions. When $\alpha=2$, the S$\alpha$S process reduces to Brownian motion, which has continuous sample paths and Gaussian increments and thus cannot capture jumps or heavy-tailed fluctuations. Accordingly, for $\alpha\in(1,2)$, S$\alpha$S processes provide a natural generalization of Brownian motion for modeling fluctuations with jumps and heavy tails.

Motivated by the aforementioned works, in this paper we study the almost sure exponential
stability of multi-dimensional stochastic systems (MSSs) driven by S$\alpha$S processes, which form a special class of L\'evy processes. The exponential stability of SDEs driven by
L\'evy processes has received considerable attention (see, e.g., Applebaum and Siakalli \cite{Applebaum2010}, Dong et al. \cite{Dong2022},  Li et al. \cite{Li2002}, Liu and Cheng \cite{Liu2023}, Zhu \cite{Zhu2018,Zhu2025}). For example, Applebaum and Siakalli \cite{Applebaum20091} used Kunita's moment estimates and exponential martingale inequality to establish criteria for almost sure exponential stability and moment exponential stability of SDEs driven by L\'evy noise. Liu et al. \cite{Liu2027} employed the generalized It\^o formula, M-matrix theory, and a contradiction argument to investigate the delay-independent mean-square exponential stability of stochastic integro-differential equations with time-varying delays and L\'evy noise.

In most of the above-mentioned literatures, the driving L\'evy process is square integrable. However, S$\alpha$S process with $1 < \alpha < 2$ is not square integrable, since its second moment is infinite. The almost sure exponential stability of multi-dimensional stochastic systems (MSSs) driven by a S$\alpha$S process has not been studied yet. The difficulty lies in the infinite variance of S$\alpha$S process, which implies that many standard tools based on second moments are no longer directly applicable, including the exponential martingale inequality under square-integrability assumptions, and quadratic Lyapunov functions.

The rest of this paper is organized as follows. In Section 2, some preliminaries on the assumptions used in this paper, $\alpha$-stable distributions, $\alpha$-stable processes, and the relationship between the fractional moments of a random variable and its characteristic function are presented. In Section 3, we establish that the solutions to MSSs almost surely remain nonzero when the initial value is nonzero. In Section 4, the $p$th moment exponential stability of solutions to MSSs is investigated. In Section 5, the almost sure exponential stability of solutions to MSSs is studied. In Section 6, a numerical example is provided to illustrate our theoretical results.

\section{Preliminaries}

Throughout this paper, unless otherwise specified, we use the following notation. Let $|\cdot|$ denote the Euclidean norm in $\mathbb{R}^n$. If $A$ is a
vector or matrix, its transpose is denoted by $A^\top$. If $A$ is a matrix,
its trace norm is denoted by $|A|=\sqrt{\text{trace}(A^\top A)}$  while its operator
norm is denoted by $||A||=\sup\{|Ax|: |x|=1\}$. Let $\mathbb{S}^{n-1}$ denote the unit sphere in $\mathbb{R}^{n}$. Let $(\Omega,\mathcal{F},\mathbb{P})$ be a complete probability space with a filtration $\{\mathcal{F}_t\}_{t\geq 0}$, which satisfies the usual conditions (i.e., it is increasing and right-continuous while $\mathcal{F}_0$ contains all $\mathbb{P}$-null sets). We write ``$\stackrel{d}{=}$'' to denote equality in distribution. Let $\{Z_t\}_{t\geq 0}$ be a standardized  one-dimensional S$\alpha$S process with $1<\alpha<2$, and its characteristic function is
\begin{align*}
\mathbb{E}[e^{iu Z_t}]=\exp\left[-t|u|^\alpha \right],
\end{align*}
and the corresponding  L\'evy measure is
\begin{align*}
\nu(dz)=\frac{C_\alpha}{|z|^{\alpha+1}}dz,
\end{align*}
where $C_\alpha=\frac{\alpha 2^{\alpha-1}}{\sqrt{\pi}}\frac{\Gamma\left(\frac{1+\alpha}{2}\right)}{\Gamma\left(1-\frac{\alpha}{2}\right)}$. Denote by $N(dt,dz)$ the Poisson random measure defined on $\mathbb{R}^+\times (\mathbb{R}-\{0\})$ with intensity measure  $\nu(dz)dt$, and let $\widetilde{N}(dt,dz)=N(dt,dz)-\nu(dz)dt$ denote the compensated Poisson random measure.

\begin{definition}\label{de2.1}
A random variable $\xi$ is said to follow a $\alpha$-stable distribution,
denoted by $\xi \sim S_\alpha (\sigma,\beta,\mu)$, if it has the characteristic function of the following
form:
 \begin{equation*}
        \Phi_\xi(u)=\mathbb{E}\exp\{iu\xi\}=\left\{ \begin{aligned}
                \exp\left\{-\sigma^\alpha|u|^\alpha\left(1-i\beta \textrm{sgn}(u)\tan\frac{\alpha\pi}{2}\right)+i\mu u\right\},\quad \textrm{if}\,\,\alpha\neq 1,\\
        \exp\left\{-\sigma |u|\left(1+i\beta\frac{2}{\pi} \textrm{sgn}(u)\log|u|\right)+i\mu u\right\},\quad \textrm{if}\,\,\alpha= 1,
                                \end{aligned} \right.
        \end{equation*}
where $\alpha\in(0,2]$, $\sigma\in(0,\infty)$, $\beta\in[-1,1]$, and $\mu\in(-\infty, \infty)$ are the index of stability, the scale, skewness, and location parameters, respectively.
\end{definition}

If $\xi \sim S_\alpha (\sigma,\beta,\mu)$ and $\alpha\in(0,2)$, then $\mathbb{E}|\xi|^p<\infty$ when $p\in(0,\alpha)$ and $\mathbb{E}|\xi|^p = \infty$ when $p\in[\alpha,\infty)$. When $\mu=0$, we say $\xi$ is strictly $\alpha$-stable. If in addition $\beta=0$, we call $\xi$ symmetric $\alpha$-stable. Note that $\xi$ is strictly 1-stable($\alpha=1$) if and only if $\beta=0$ (symmetric case). We refer to  Janicki and Weron \cite{Janicki}, and Sato  \cite{Sato} for more details on $\alpha$-stable distributions.

\begin{definition}
An $\mathcal{F}_t$-adapted stochastic process $\{L_t\}_{t\ge0}$ is called a standard $\alpha$-stable process if \\
(1) $L_0 = 0$, a.s. \\
(2) $L_t - L_s \sim S_\alpha ((t-s)^{\frac{1}{\alpha}}, \beta, 0)$, $t> s \ge 0$; \\
(3) For any finite tume points $0\le s_0 < s_1 < ...<s_m <\infty$, the random variables
$L_{s_0}, L_{s_1}-L_{s_0}, ... , L_{s_m}-L_{s_{m-1}}$ are independent.
\end{definition}

Note that a standard  $\alpha$-stable process has the self-similarity property $L_{at}\stackrel{d}{=}a^{\frac{1}{\alpha}}L_t$ with $a>0$. Every L\'evy process has a unique modification which is c\`{a}dl\`{a}g (right continuous with left limits) and which is also a L\'evy process. We will henceforth assume that we are using the c\`{a}dl\`{a}g version of the given $\alpha$-stable process.

Consider the following MSS
\begin{align}\label{model}
dy(t) = f(y(t))dt + g(y(t))dZ_t,\qquad  y(0)=y_0\in\mathbb{R}^n,
\end{align}
where $f,g: \mathbb{R}^n\rightarrow \mathbb{R}^n$. We impose the following assumptions.

\begin{assumption}\label{ass3.1}
Assume that $f,g: \mathbb{R}^n\rightarrow \mathbb{R}^n$ are globally Lipschitz continuous, that is, there exist constants $L_f, L_g>0$ such that
\begin{align*}
|f(x)-f(y)|\leq L_f|x-y|,\quad |g(x)-g(y)|\leq L_g|x-y|\quad x,y\in \mathbb{R}^n.
\end{align*}
Moreover,  $f(0)=0$ and $g(0)=0$, which guarantee that $y(t)\equiv0$ is the trivial solution corresponding to the initial value $y_0=0$.
\end{assumption}

Under Assumption \ref{ass3.1}, there exists a unique c\`{a}dl\`{a}g adapted strong solution $y(t)$ for every $y_0\in\mathbb{R}^n$.

\begin{assumption}\label{ass3.2}
There exist constants $0<\rho_1\leq L_g$ satisfying
\begin{align*}
 |y^Tg(y)|^2\geq \rho_1|y|^4,\quad \forall y\in \mathbb{R}^n.
\end{align*}
In particular, $g(y)\neq 0$ for every $y\neq 0$.
\end{assumption}

The following lemma, which plays an important role in the proofs of our main results, states that the moments of a random variable can be obtained by using the relationship between its fractional moments and its characteristic function.

\begin{lemma}(see, Proposition 6,  Fageot et al. \cite{Fageot2017})\label{le2.1}
For a random variable $X$ with characteristic function $\Phi_X$ and $s\in(0,2)$, we have the relation
\begin{align*}
\mathbb{E}|X|^s=\left(\int_\mathbb{R}\frac{1-\cos \nu}{|\nu|^{s+1}}d\nu\right)^{-1}\int_\mathbb{R}\frac{1-\Re(\Phi_X)(\nu)}{|\nu|^{s+1}}d\nu\in[0,\infty],
\end{align*}
where $\mathcal{R}(z)$ denotes the real part of $z\in \mathbb{C}$.
\end{lemma}

\section{Non-vanishing of solutions to MSS (\ref{model})}

In this section, we show that the solution to the MSS (\ref{model}) with a nonzero initial value almost surely remains nonzero under Assumption \ref{ass3.1}. To this end, we first introduce the following elementary inequality.

\begin{lemma}\label{le3.1}
Let $p>0$ and $M>0$. Then, there exists a constant $K=K(p,M)>0$ such that for all $u\in\mathbb{S}^{n-1}$ and $v\in \mathbb{R}^n$ with $|v|\leq M$, and $u+v\neq0$,
\begin{align}\label{3.1}
\bigl||u+v|^{-p}-1+pu^\top v\bigr|\leq K\frac{|v|^2}{|u+v|^{p}},
\end{align}
where
\begin{align*}
    K=\max\left\{(1+M)^p\frac{p(p+3)}2 2^{p+2}, 4\left[1+(1+pM)(1+M)^p\right] \right\}.
\end{align*}
\end{lemma}

\begin{proof}
If $v=0$, the inequality is trivial. Hence we only need to consider $v\neq0$.
Next, we prove the (\ref{3.1}) by considering two cases.

\textbf{Case 1:} $0<|v|\le\frac12$. Define $d(\vartheta)=|u+\vartheta v|^{-p}$, $\vartheta\in[0,1]$.
Then $d(0)=1$, $d'(0)=-p u^\top v$, and
\begin{align*}
|d''(\vartheta)|
= \bigl| p(p+2)|u+\vartheta v|^{-p-4}((u+\vartheta v)\cdot v)^2 - p|u+\vartheta v|^{-p-2}|v|^2 \bigr|
\le p(p+3)2^{p+2}|v|^2.
\end{align*}
By Taylor's formula, there exists $\varsigma\in(0,1)$ such that $d(1)=d(0)+d'(0)+\frac12 d''(\varsigma),$
which implies
\begin{align*}
\bigl||u+v|^{-p}-1+p u^\top v\bigr| \le \frac{p(p+3)}2 2^{p+2}|v|^2.
\end{align*}
Since $|u+v|\le|u|+|v|\le1+M$, we have
\begin{align}\label{case1}
\frac{|u+v|^p}{|v|^2}\bigl||u+v|^{-p}-1+p u^\top v\bigr|
\le (1+M)^p\frac{p(p+3)}2 2^{p+2}.
\end{align}

\noindent\textbf{Case 2:} $|v|>\frac12$. Then,
\begin{align}\label{case2}
\frac{|u+v|^p}{|v|^2}
\left||u+v|^{-p}-1+pu^\top v\right|&=\frac{1}{|v|^2}
\left|1-|u+v|^p+pu^\top v\,|u+v|^p\right|\nonumber\\
&\leq4\left[1+(1+pM)(1+M)^p\right].
\end{align}
Combining (\ref{case1}) and (\ref{case2}) yields (\ref{3.1}).
This completes the proof.
\end{proof}

\begin{theorem}\label{th3.1}
If the initial value satisfies $y_0\neq 0$ and Assumptions \ref{ass3.1} and \ref{ass3.2} hold, then
\begin{align}\label{3.3}
\mathbb{P}(y(t)\neq 0\,\, for \,\, all \,\, t\geq 0)=1.
\end{align}
\end{theorem}

\begin{proof}
Assume that (\ref{3.3}) is false. This implies that, for some $y_0\neq 0$, there will be a stopping time $\iota$ with $\mathbb{P}(\iota<\infty)>0$ when the solution will be zero for the first time:
\begin{align*}
\iota=\inf\{t\geq 0: |y(t)|=0\}.
\end{align*}
Since the paths of $y$ are almost-surely right continuous with left limits, we find a pair of constants  $T>0$ and $\kappa>1$ sufficiently large for $\mathbb{P}(\Lambda)>0$, where
\begin{align*}
\Lambda=\left\{\omega\in\Omega: \iota(\omega)\leq T\quad \mathrm{and} \quad |y(t,\omega)|\leq \kappa -1\quad \mathrm{for} \quad \mathrm{all}\quad 0\leq t\leq \iota(\omega) \right\}.
\end{align*}
Let $V(y)=|y|^{-p}$ for $y\neq 0$, where $p\in(0,\alpha-1)$. Applying It\^o's formula, we obtain
\begin{align}\label{3.6}
dV(y(t)) = \mathcal{L}V(y(t))dt + \int_{\mathbb{R}} \bigl( |y(t^-)+g(y(t^-))z|^{-p} - |y(t^-)|^{-p} \bigr)\widetilde{N}(dz,dt),
\end{align}
where
\begin{align*}
\mathcal{L}V(y) &= -p|y|^{-p-2} y^\top f(y)+ \int_{\mathbb{R}} \Bigl( |y+g(y)z|^{-p} - |y|^{-p} + p|y|^{-p-2} y^\top g(y)z \Bigr)\nu(dz).
\end{align*}
It follows from Assumption \ref{ass3.1} that
\begin{align}\label{3.7}
\mathcal{L}V(y)&=-p|y|^{-p-2} y^\top f(y)+\frac{C_\alpha}{2}\int_{\mathbb{R}}\left(|y+g(y)z|^{-p}+|y-g(y)z|^{-p}-2|y|^{-p}\right)\frac{1}{|z|^{\alpha+1}}\nonumber\\
&\leq |y|^{-p}\left(pL_f+\frac{C_\alpha}{2}\int_{\mathbb{R}}\left||u+\varpi z|^{-p}+|u-\varpi z|^{-p}-2|\right|\frac{1}{|z|^{\alpha+1}}dz\right),
\end{align}
where $u=y/|y|\in \mathbb{S}^{n-1}$, and $\varpi=\frac{g(y)}{|y|}$. Let
\begin{align}\label{omega}
\Omega_1&=\{|z|\leq \delta\}, \quad\Omega_2=\{\delta<|z|\leq R_0\}\setminus \Omega_3,\nonumber\\
\Omega_3&=\left[-\frac{1}{\sqrt{\rho_1}}-\varepsilon_0,\frac{1}{\sqrt{\rho_1}}+\varepsilon_0\right]\cap\{\delta< |z|\leq R_0\},\quad \Omega_4=\{ |z|>R_0\}.
\end{align}
where $\delta=\frac{1}{2\sqrt{L_g}}$, $R_0=\max\{\frac{4}{\sqrt{\rho_1}},\delta+1\}$, and $\varepsilon_0=\frac{1}{4}\min\{\delta,1/\sqrt{\rho_1}\}$. Obviously, $\mathbb{R}=\bigcup_{i=1}^4\Omega_i$. By  Cauchy-Schwarz inequality and Assumption \ref{ass3.2}, we have
\begin{align}\label{r3.1}
 \rho_1\leq \mu_1^2\leq \mu_2\leq L_g,
\end{align}
where $\mu_1=u^\top\varpi$ and $\mu_2=|\varpi|^2$. Note that $|u+\varpi z|\geq |u|-|\varpi z|\geq \frac{1}{2}$ on $\Omega_1$. Applying Lemma \ref{le3.1} with $M=\frac{1}{2}$ and using (\ref{r3.1}), we obtain
\begin{align}\label{3.9}
\int_{\Omega_1}\left||u+\varpi z|^{-p}+|u-\varpi z|^{-p}-2|\right|\frac{1}{|z|^{\alpha+1}}dz&\leq
\int_{\Omega_1}\left||u+\varpi z|^{-p}-1+p\mu_1z\right|\frac{1}{|z|^{\alpha+1}}dz\nonumber\\
&\quad +
\int_{\Omega_1}\left||u-\varpi z|^{-p}-1-p\mu_1z\right|\frac{1}{|z|^{\alpha+1}}dz\nonumber\\
&\leq K\int_{\Omega_1}\left(\frac{|z\varpi|^2}{|u+z\varpi|^p}+\frac{|z\varpi|^2}{|u-z\varpi|^p}\right)\frac{1}{|z|^{\alpha+1}}dz\nonumber\\
&\leq 2^{p+1} K L_g\int_{\Omega_1}\frac{1}{|z|^{\alpha-1}}dz\nonumber\\
&= \frac{2^{p+2} K L_g}{2-\alpha}\delta^{2-\alpha}.
\end{align}
For $z\in \Omega_2\bigcup\Omega_3$, by (\ref{r3.1}), we have
\begin{align}\label{o2}
|u\pm z\varpi|^2=\mu_2\left(z\pm\frac{\mu_1}{\mu_2}\right)^2+1-\frac{\mu_1^2}{\mu_2}\geq \rho_1\left( z\pm\frac{\mu_1}{\mu_2}\right)^2.
\end{align}
Then,
\begin{align}\label{3.10}
\int_{\Omega_2}\left||u+\varpi z|^{-p}+|u-\varpi z|^{-p}-2|\right|\frac{1}{|z|^{\alpha+1}}dz&\leq 2\left(\rho_1^{-\frac{p}{2}}\varepsilon_0^{-p}+1\right)\int_{\frac{1}{\sqrt{\rho_1}}+\varepsilon_0\leq |z|\leq R_0}\frac{1}{|z|^{\alpha+1}}dz\nonumber\\
& \leq \frac{4\left(\rho_1^{-\frac{p}{2}}\varepsilon_0^{-p}+1\right)}{\alpha}\left[\left(\frac{1}{\sqrt{\rho_1}}+\varepsilon_0\right)^{-\alpha}-R_0^{-\alpha}\right].
\end{align}
and
\begin{align}\label{3.11}
&\int_{\Omega_3}\left||u+\varpi z|^{-p}+|u-\varpi z|^{-p}-2|\right|\frac{1}{|z|^{\alpha+1}}dz\nonumber\\
&\leq \delta^{-\alpha-1}\rho_1^{-p/2}\int_{\Omega_3}\left(\left| z+\frac{\mu_1}{\mu_2}\right|^{-p}+\left| z-\frac{\mu_1}{\mu_2}\right|^{-p}\right)dz+2 \int_{\Omega_3}\frac{1}{|z|^{\alpha+1}}dz\nonumber\\
&\leq2\delta^{-\alpha-1}\rho_1^{-p/2}\int_{-\frac{2}{\sqrt{\rho_1}}-\varepsilon_0}^{\frac{2}{\sqrt{\rho_1}}+\varepsilon_0}|s|^{-p}ds+2 \int_{\delta<|z|\leq R_0}\frac{1}{|z|^{\alpha+1}}dz \nonumber\\
&\leq\frac{4 \delta^{-\alpha-1}\rho_1^{-p/2}\left(\frac{2}{\sqrt{\rho_1}}+\varepsilon_0\right)^{1-p}}{1-p}+\frac{4}{\alpha}\delta^{-\alpha}.
\end{align}
For $z\in \Omega_4$, it follows from (\ref{r3.1}) that
\begin{align*}
|u\pm z\varpi|\geq |z||\varpi|-1\geq \sqrt{\rho_1}|z|-1\geq \frac{\sqrt{\rho_1}}{2}|z|.
\end{align*}
This implies that
\begin{align}\label{3.12}
\int_{\Omega_4}\left||u+\varpi z|^{-p}+|u-\varpi z|^{-p}-2|\right|\frac{1}{|z|^{\alpha+1}}dz&\leq 2\int_{\Omega_4}\left(2^{p}\rho_1^{-p/2}|z|^{-p}+1\right)\frac{1}{|z|^{\alpha+1}}dz\nonumber\\
&\leq4\int_{R_0}^\infty\left(2^{p}\rho_1^{-p/2}r^{-p-\alpha-1}+r^{-\alpha-1}\right)dr\nonumber\\
&=\frac{2^{p+2}\rho_1^{-\frac{p}{2}}}{p+\alpha}R_0^{-p-\alpha}+\frac{4}{\alpha}R_0^{-\alpha}.
\end{align}
By (\ref{3.7}), (\ref{3.9}), and (\ref{3.10})-(\ref{3.12}), we can conclude that
\begin{align*}
\mathcal{L}V(y)\leq \left(pL_f+\frac{C_\alpha}{2}K_1\right)V(y),
\end{align*}
where
\begin{align*}
K_1&=2^{p+2} K L_g\frac{\delta^{2-\alpha}}{2-\alpha}+\frac{4\left(\rho_1^{-\frac{p}{2}}\varepsilon_0^{-p}+1\right)}{\alpha}\left[\left(\frac{1}{\sqrt{\rho_1}}+\varepsilon_0\right)^{-\alpha}-R_0^{-\alpha}\right]\nonumber\\
&\quad+\frac{4 \delta^{-\alpha-1}\rho_1^{-p/2}\left(\frac{2}{\sqrt{\rho_1}}+\varepsilon_0\right)^{1-p}}{1-p}+\frac{4}{\alpha}\delta^{-\alpha}+\frac{2^{p+2}\rho_1^{-\frac{p}{2}}}{p+\alpha}R_0^{-p-\alpha}+\frac{4}{\alpha}R_0^{-\alpha}.
\end{align*}
Define
\begin{align*}
\iota_\varepsilon=\inf\{t\geq 0: |y(t)|\leq \varepsilon\quad \mathrm{or}\quad |y(t)|\geq \kappa\}
\end{align*}
for each $0<\varepsilon<|y_0|$. Then, by It\^o's formula,
\begin{align*}
&\mathbb{E}[e^{-(pL_f+C_\alpha K_1)(\iota_\varepsilon\wedge T)}V(y(\iota_\varepsilon\wedge T))]\nonumber\\
&=V(y_0)+\mathbb{E}\left[ \int_0^{\iota_\varepsilon\wedge T}e^{-(pL_f+C_\alpha K_1) s}\left(-(pL_f+C_\alpha K_1)V(y(s))+\mathcal{L}V(y(s))\right)ds\right]\nonumber\\
&\leq |y_0|^{-p}.
\end{align*}
Note that, if $\omega\in \Lambda$ then $\iota_\varepsilon(\omega)\leq T$ and $|y(\iota_\varepsilon(\omega))|\leq \varepsilon$. Then,
\begin{align*}
e^{-(pL_f+C_\alpha K_1) T}\varepsilon^{-p}\mathbb{P}(\Lambda)&\leq \mathbb{E}[e^{-(pL_f+C_\alpha K_1) T}\varepsilon^{-p}I_{\Lambda} ]\nonumber\\
&\leq \mathbb{E}[e^{-(pL_f+C_\alpha K_1)\iota_\varepsilon}|y(\iota_\varepsilon)|^{-p}I_{\Lambda} ]\nonumber\\
&=\mathbb{E}[e^{-(pL_f+C_\alpha K_1)(\iota_\varepsilon\wedge T)}V(y(\iota_\varepsilon\wedge T))I_{\Lambda} ]\nonumber\\
&\leq \mathbb{E}[e^{-(pL_f+C_\alpha K_1)(\iota_\varepsilon\wedge T)}V(y(\iota_\varepsilon\wedge T))]\nonumber\\
&\leq |y_0|^{-p}.
\end{align*}
Hence,
\begin{align}
\mathbb{P}(\Lambda)\leq \varepsilon^p e^{(pL_f+C_\alpha K_1)T}|y_0|^{-p}.
\end{align}
Since $0<\varepsilon<|y_0|$ and $\varepsilon>0$ is arbitrary, we must have $\mathbb{P}(\Lambda)=0$, but this contradicts the definition of $\Lambda$. The proof is complete.
\end{proof}

\section{$p$th moment exponential stability of solutions to MSS (\ref{model})}
In this section, we establish the $p$th moment exponential stability for the MSS (\ref{model}). To this end, we introduce the following two key lemmas.

\begin{lemma} \label{le2.3}For any $s\in(1,2)$, we have
\begin{align}\label{2.3main}
\int_0^1 u^{s-1}\cos u\,du - \sin 1 - (s-1)\cos 1 - (s-1)(s-2)\int_1^\infty u^{s-3}\cos u\,du = \Gamma(s)\cos\frac{\pi s}{2}.
\end{align}
\end{lemma}

\begin{proof} Let $\varepsilon>0$ and $I(\varepsilon)=\int_0^\infty e^{-\varepsilon u}u^{s-1}\cos udu $. Since $\lim\limits_{u\rightarrow\infty}u^{s-1}e^{-\frac{\varepsilon}{2}u}=0$, there exists a constant $\bar{C}>0$ such that $u^{s-1}e^{-\frac{\varepsilon}{2}u}\leq \bar{C}$ for all $u\geq 1$. It follows that
\begin{align}\label{c1.1}
\left|\int_1^\infty e^{-\varepsilon u}u^{s-1}\cos udu\right|&\leq \int_1^\infty e^{-\varepsilon u}u^{s-1}du\nonumber\\
&\leq \bar{C}\int_1^\infty e^{-\frac{\varepsilon}{2} u}du\nonumber\\
&=\bar{C}\frac{2}{\varepsilon}e^{-\frac{\varepsilon}{2}}.
\end{align}
By (\ref{c1.1}) and
\begin{align}\label{c1.2}
\left|\int_0^1 e^{-\varepsilon u}u^{s-1}\cos udu\right|&\leq \int_0^1 u^{s-1}du=\frac{1}{s},
\end{align}
 we get that $I(\varepsilon)$ is absolutely convergent. Combining $\cos u = \Re(e^{iu})$ with the Laplace-Gamma integral formula gives that
\begin{align}\label{c1.3}
I(\varepsilon)=\Re\int_0^\infty e^{-(\varepsilon-i)u}u^{s-1}du=\Re\bigl(\Gamma(s)(\varepsilon-i)^{-s}\bigr).
\end{align}
By  $\ln(-i)=-\frac{\pi i}{2}$, Euler's formula, and (\ref{c1.3}), we have
\begin{align}\label{I2.0}
\lim\limits_{\varepsilon\rightarrow 0^+}I(\varepsilon)=\Re\bigl(\Gamma(s)(-i)^{-s}\bigr)=\Re\bigl(\Gamma(s)e^{-s\ln(-i)}\bigr)=\Re\bigl(\Gamma(s)e^{i\frac{\pi s}{2}}\bigr)=\Gamma(s)\cos\frac{\pi s}{2}.
\end{align}
Note that
\begin{align}\label{I2.1}
I(\varepsilon)&=\int_0^1 e^{-\varepsilon u}u^{s-1}\cos u\,du+\int_1^\infty e^{-\varepsilon u}u^{s-1}\cos u\,du\nonumber\\
&:=I_1(\varepsilon)+I_2(\varepsilon).
\end{align}
For $u\in[0,1]$ and any $\varepsilon>0$, we have $|e^{-\varepsilon u}u^{s-1}\cos u|\le u^{s-1}$, and $\int_0^1 u^{s-1}\,du = \frac{1}{s} < \infty$.
The dominated convergence theorem implies that
\begin{align}\label{I2.2}
\lim_{\varepsilon\to0^+} I_1(\varepsilon)
= \int_0^1 \lim_{\varepsilon\to0^+} e^{-\varepsilon u}u^{s-1}\cos u \, du
= \int_0^1 u^{s-1}\cos u \, du .
\end{align}
Now we consider $I_2(\varepsilon)$. By integrating by parts, it follows that
\begin{align}\label{I2.3}
I_2(\varepsilon) &= e^{-\varepsilon u}u^{s-1}\sin u\biggl|_1^\infty
-\int_1^\infty \sin u\cdot\frac{d}{du}\!\big(e^{-\varepsilon u}u^{s-1}\big)du \nonumber\\
&= -e^{-\varepsilon}\sin 1+\varepsilon\int_1^\infty e^{-\varepsilon u}u^{s-1}\sin u\,du
-(s-1)\int_1^\infty e^{-\varepsilon u}u^{s-2}\sin u\,du .
\end{align}
A further integration by parts gives
\begin{align}\label{I2.4}
\int_1^\infty e^{-\varepsilon u}u^{s-1}\sin udu&=-e^{-\varepsilon u}u^{s-1}\cos u\biggl|_{1}^\infty+\int_1^\infty\cos ud(e^{-\varepsilon u}u^{s-1})\nonumber\\
&=e^{-\varepsilon}\cos 1-\varepsilon\int_1^\infty e^{-\varepsilon u}u^{s-1}\cos udu+(s-1)\int_1^\infty e^{-\varepsilon u}u^{s-2}\cos udu
\end{align}
and
\begin{align}\label{I2.5}
\int_1^\infty e^{-\varepsilon u}u^{s-2}\sin udu
&= -e^{-\varepsilon u}u^{s-2}\cos u\biggl|_1^\infty
+\int_1^\infty \cos u\cdot\frac{d}{du}\!\big(e^{-\varepsilon u}u^{s-2}\big)du \nonumber\\
&= e^{-\varepsilon}\cos 1 - \varepsilon\int_1^\infty e^{-\varepsilon u}u^{s-2}\cos udu + (s-2)\int_1^\infty e^{-\varepsilon u}u^{s-3}\cos udu .
\end{align}
By combining (\ref{I2.3})-(\ref{I2.5}), we can conclude that
\begin{align}\label{I2.6}
I_2(\varepsilon) &= -e^{-\varepsilon}\sin 1 - (s-1)e^{-\varepsilon}\cos 1 + \varepsilon e^{-\varepsilon}\cos 1-\varepsilon^2\int_1^\infty e^{-\varepsilon u}u^{s-1}\cos udu \nonumber\\
&\quad + 2(s-1)\varepsilon\int_1^\infty e^{-\varepsilon u}u^{s-2}\cos udu - (s-1)(s-2)\int_1^\infty e^{-\varepsilon u}u^{s-3}\cos udu.
\end{align}
For $u\in[1,\infty)$ and any $\varepsilon>0$, we have
$|e^{-\varepsilon u}u^{s-3}\cos u| \le u^{s-3}$,
and $\int_1^\infty u^{s-3}\,du = \frac{1}{2-s} < \infty$ (since $s<2$).
By the dominated convergence theorem,
\begin{align}\label{I2.7}
\lim_{\varepsilon\to0^+} \int_1^\infty e^{-\varepsilon u}u^{s-3}\cos u\,du
= \int_1^\infty \lim_{\varepsilon\to0^+} \bigl(e^{-\varepsilon u}u^{s-3}\cos u\bigr)\,du
= \int_1^\infty u^{s-3}\cos u\,du .
\end{align}
Since $s\in(1,2)$, we have
\begin{align}
\left|\varepsilon^2\int_1^\infty e^{-\varepsilon u}u^{s-1}\cos udu\right|&\leq \varepsilon^2\int_0^\infty e^{-\varepsilon u}u^{s-1}du=\Gamma(s)\varepsilon^{2-s}\label{I2.8},\\
\left|\varepsilon\int_1^\infty e^{-\varepsilon u}u^{s-2}\cos udu\right|&\leq \varepsilon \int_0^\infty e^{-\varepsilon u}u^{s-2}du=\Gamma(s-1)\varepsilon^{2-s}\label{I2.9},
\end{align}
both of which tend to zero as $\varepsilon\rightarrow 0^+$. By (\ref{I2.6})-(\ref{I2.9}), we get that
\begin{align}\label{I2.10}
\lim_{\varepsilon\to0^+}I_2(\varepsilon)=-\sin 1 - (s-1)\cos 1 - (s-1)(s-2)\int_1^\infty u^{s-3}\cos udu.
\end{align}
It follows from (\ref{I2.1}), (\ref{I2.2}), and (\ref{I2.10}) that
\begin{align}\label{I2.11}
\lim_{\varepsilon\to0^+}I(\varepsilon)=\int_0^1 u^{s-1}\cos u \, du-\sin 1 - (s-1)\cos 1 - (s-1)(s-2)\int_1^\infty u^{s-3}\cos udu.
\end{align}
Combining (\ref{I2.0}) with (\ref{I2.11}) gives that (\ref{2.3main}) holds. This completes the proof.
\end{proof}

\begin{lemma} \label{le2.4}
Assume that $p\in(0,\alpha-1)$, $\alpha\in(1,2)$, and $\vartheta\neq 0$. Let
\begin{align}\label{gt}
\varphi(t)=\frac{2\Gamma(1+p)\sin\frac{\pi p}{2}}{\pi}\int_0^\infty\frac{1-e^{-|\vartheta|^\alpha tr^\alpha}\cos r}{r^{p+1}}dr.
\end{align}
Then, $\varphi(0)=1$, $\lim\limits_{t\rightarrow \infty}\varphi(t)=\infty$, $\varphi'(0^+)=-p\widehat{H}_\alpha(p)|\vartheta|^\alpha$, where
\begin{align}\label{Hal}
\widehat{H}_\alpha(p)=-\frac{2\Gamma(1+p)\sin\frac{p\pi}{2}}{\pi p}\Gamma(\alpha-p)\cos\frac{(\alpha-p)\pi}{2}
\end{align}
\end{lemma}

\begin{proof}
 Note that
\begin{align}\label{G}
\int_0^\infty\frac{1-\cos u}{u^{p+1}}du=\frac{\pi}{2\Gamma(p+1)\sin\frac{\pi p}{2}}.
\end{align}
Then, by the definition of $\varphi$, we have $\varphi(0)=1$.  Let $t_0\geq 0$ be arbitrary and sequence $\{t_m\}_{m=1}^\infty\subset[0,\infty)$ with $\lim\limits_{m\rightarrow\infty}t_m=t_0$. This implies that there exists $T>0$ such that $t_m,t_0\in[0,T]$ for all $m$. Let $q(t,u)=\frac{1-e^{-|\vartheta|^\alpha tu^\alpha}\cos u}{u^{p+1}}$. It is easy to see that $q(t,u)\geq 0$ and $\lim_{m\rightarrow\infty}q(t_m,u)=q(t_0,u)$. Applying
\begin{align*}
1-\cos u\leq \frac{u^2}{2},\quad 1-e^{-u}\leq u, \quad u\geq 0,
\end{align*}
we get that for all $t_m\in[0,T]$,
\[
q(t_m,u)\leq h(u):=
\begin{cases}
\frac{1}{2}u^{1-p}+|\vartheta|^\alpha Tu^{\alpha-p-1}, & 0 < u \leq 1, \\
\frac{2}{u^{p+1}}, & u > 1.
\end{cases}
\]
By dominated convergence theorem, and
\begin{align*}
\int_0^\infty h(u)du=\frac{1}{2(2-p)}+\frac{|\vartheta|^\alpha T}{\alpha-p}+\frac{2}{p},
\end{align*}
we have $\lim_{t\rightarrow t_0}\varphi(t)=\varphi(t_0)$. Since $t_0$ is arbitrary, $\varphi(t)$ is continuous on $[0,\infty)$. For any sequence $\{t^{(1)}_m\}_{m=1}^\infty\subset[0,\infty)$ with $\lim\limits_{m\rightarrow\infty}t^{(1)}_m=\infty$, by Fatou's lemma,
\begin{align*}
\liminf\limits_{m\rightarrow\infty}\varphi(t^{(1)}_m)&=\frac{2\Gamma(1+p)\sin\frac{\pi p}{2}}{\pi}\liminf\limits_{m\rightarrow\infty}\int_0^\infty q(t^{(1)}_m,u)du\\
&\geq \frac{2\Gamma(1+p)\sin\frac{\pi p}{2}}{\pi}\int_0^\infty \liminf\limits_{m\rightarrow\infty}q(t^{(1)}_m,u)du\\
&=\frac{2\Gamma(1+p)\sin\frac{\pi p}{2}}{\pi}\int_0^\infty\frac{1}{u^{p+1}}du=\infty.
\end{align*}
Then, $\lim\limits_{t\rightarrow\infty}\varphi(t)=\infty$. It follows from the definition of $q(t,u)$ that
\begin{align*}
\partial_t q(t,u) = e^{-|\vartheta|^\alpha t u^\alpha} |\vartheta|^\alpha u^{\alpha-p-1} \cos u.
\end{align*}
Fix an arbitrary $t_1>0$. For every $t\in[\frac{t_1}{2}, \frac{3t_1}{2}]$, we have $|\partial_t q(t,u)| \leq |\vartheta|^\alpha e^{-|\vartheta|^\alpha \frac{t_1}{2} u^\alpha} u^{\alpha-1-p}$, and
\begin{align*}
\int_0^\infty |\vartheta|^\alpha e^{-|\vartheta|^\alpha \frac{t_1}{2} u^\alpha} u^{\alpha-1-p}du=\frac{|\vartheta|^p}{\alpha}t_1^{\frac{p}{\alpha}-1}2^{1-\frac{p}{\alpha}}\Gamma\left(1-\frac{p}{\alpha}\right).
\end{align*}
Hence,  $\varphi(t)$ is differentiable at $t_1$. Since $t_1>0$ is arbitrary, $\varphi(t)$ is differential on $(0,\infty)$ and
 \begin{align}\label{u0}
\varphi'(t)&=\frac{2\Gamma(1+p)\sin\frac{\pi p}{2}}{\pi}|\vartheta|^\alpha \int_0^\infty G(u,t) \cos udu,\quad t> 0,
\end{align}
where $G(u,t)=e^{-|\vartheta|^\alpha t u^\alpha} u^{\alpha-p-1}$. It is obvious to see that
\begin{align}\label{u1}
|G(u,t)\cos u|\leq u^{\alpha-p-1},\quad u\in(0, 1).
\end{align}
Since $\int_0^1 u^{\alpha-p-1}du=\frac{1}{\alpha-p}$, by (\ref{u1}) and dominated convergence theorem, it follows that
\begin{align}\label{u2}
\lim\limits_{t\rightarrow 0}\int_0^1G(u,t)\cos udu&=\int_0^1u^{\alpha-p-1}\cos udu.
\end{align}
Integrating by parts twice, we obtain
\begin{align}\label{u4}
\int_1^\infty G(u,t) \cos  u du&=-G(1,t)\sin 1-\int_1^\infty\sin u \frac{\partial G(u,t)}{\partial  u}du\nonumber\\
&=-G(1,t)\sin 1-\cos 1G_u(1,t) -\int_1^\infty \cos u\frac{\partial^2 G(u,t)}{\partial  u^2}du,
\end{align}
where
\begin{align*}
G(1,t)&=e^{-|\vartheta|^\alpha t},\\
G_u(1,t)&=e^{-|\vartheta|^\alpha t}(\alpha-p-1-\alpha|\vartheta|^\alpha t),\\
\frac{\partial^2 G(u,t)}{\partial  u^2}&=e^{-|\vartheta|^\alpha t u^\alpha} u^{\alpha-p-3}\left[(\alpha-p-1)(\alpha-p-2)-\alpha |\vartheta|^\alpha (3\alpha-2p-3)u^\alpha t+\alpha^2|\vartheta|^{2\alpha}t^2 u^{2\alpha}\right].
\end{align*}
Applying $\sup_{y\geq 0}y^ke^{-y}=\left(\frac{k}{e}\right)^{k}$, $k=0,1,2$, we get for any $t>0$ and $u\geq 1$,
\begin{align}\label{u3}
\left|\frac{\partial^2 G(u,t)}{\partial  u^2}\cos u \right|&\leq e^{-|\vartheta|^\alpha t u^\alpha}u^{\alpha-p-3}\left(|(\alpha-p-1)(\alpha-p-2)|+\alpha|3\alpha-2p-3||\vartheta|^\alpha t u^\alpha\right.\nonumber\\
&\left.\quad+\alpha^2 |\vartheta|^{2\alpha} t^2 u^{2\alpha}\right)\nonumber\\
&\leq \left(|(\alpha-p-1)(\alpha-p-2)|+\frac{\alpha|3\alpha-2p-3|}{e}+\frac{4\alpha^2}{e^2}\right)u^{\alpha-p-3}.
\end{align}
Since $\int_1^\infty u^{\alpha-p-3}du=\frac{1}{2-\alpha+p}$, combining (\ref{u3}) with dominated convergence theorem yields that
\begin{align}\label{u5}
\lim\limits_{t\rightarrow 0^+}\int_1^\infty \frac{\partial^2 G(u,t)}{\partial  u^2}\cos udu&=(\alpha-1-p)(\alpha-p-2)\int_1^\infty u^{\alpha-p-3}\cos udu.
\end{align}
By (\ref{u2}), (\ref{u4}), and (\ref{u5}), we can conclude that
\begin{align}\label{u6}
\lim\limits_{t\rightarrow 0^+}\int_0^\infty G(u,t) \cos udu&=\lim\limits_{t\rightarrow 0^+}\int_0^1 G(u,t) \cos udu+\lim\limits_{t\rightarrow 0^+}\int_1^\infty G(u,t) \cos udu\nonumber\\
&=\int_0^1u^{\alpha-p-1}\cos udu-\sin1-(\alpha-p-1)\cos 1\nonumber\\
&\quad-(\alpha-1-p)(\alpha-p-2)\int_1^\infty u^{\alpha-p-3}\cos udu.
\end{align}
Since $p\in(0,\alpha-1)$ and $\alpha\in(1,2)$,
we have $\alpha-p\in(1,2)$. By (\ref{u0}), (\ref{u6}), and Lemma \ref{le2.3}, we find
\begin{align}\label{u10}
\lim\limits_{t\rightarrow 0^+}\varphi'(t)=\frac{2\Gamma(1+p)\sin\frac{\pi p}{2}}{\pi}|\vartheta|^\alpha\Gamma(\alpha-p)\cos\frac{(\alpha-p)\pi}{2}.
\end{align}
For any $t>0$, using the mean value theorem,
 there exists $t_\xi\in(0,t)$ such that $\varphi'(t_\xi)=\frac{\varphi(t)-\varphi(0)}{t}$. Letting $t\rightarrow 0^+$ gives $t_\xi\rightarrow 0^+$, and (\ref{u10}) yields
\begin{align}\label{g0}
\varphi'(0^+)=\lim\limits_{t\rightarrow 0^+}\frac{\varphi(t)-\varphi(0)}{t}=-p\widehat{H}_\alpha(p)|\vartheta|^\alpha,
\end{align}
where $\widehat{H}_\alpha(p)$ is given in (\ref{Hal}). Since $0<p<\alpha-1<1$ and $1<\alpha-p<2$,we have $\widehat H_\alpha(p)>0$. This completes the proof.
\end{proof}

We now proceed to establish the main result of this section.

\begin{theorem}\label{th4.1}
Suppose that Assumptions \ref{ass3.1} and \ref{ass3.2} hold. If
\begin{align}\label{gamm}
-\Gamma(\alpha)\cos\frac{\alpha\pi}{2}\rho_1^{\frac{\alpha}{2}}-L_f-K_2\sqrt{L_g-\rho_1}-K_3(L_g-\rho_1)>0,
\end{align}
where
\begin{align}\label{K2}
K_2&=C_\alpha\left(\frac{4}{3}\right)^{1+\alpha}\left[\frac{25}{12(2-\alpha)}\left(\max\left\{1,L_g^{\frac{\alpha-1}{2}}\right\}\right)^2\right.\nonumber\\
&\quad\left.+4\max\left\{1,\left(\frac{4\sqrt{L_g}}{\sqrt{\rho_1}}\right)^{\alpha-1}\right\}\left(4+\frac{5}{2}\log\frac{5}{4}\right)L_g^{\frac{\alpha-1}{2}}\right]
\end{align}
and
\begin{align}\label{K3}
K_3&=\frac{16C_\alpha}{9(2-\alpha)}4^{\alpha-2}L_g^{\frac{\alpha-2}{2}}+\frac{256C_\alpha}{\alpha\rho_1^2}\max\{L_g,\rho_1\}\left(4^\alpha L_g^{\frac{\alpha}{2}}-4^{-\alpha}\rho_1^{\frac{\alpha}{2}}\right)\nonumber\\
&\quad+C_\alpha\max\left\{\frac{4}{\rho_1},\frac{2^{3-\alpha}}{\rho_1^{\frac{3-\alpha}{2}}}\right\} 4^{-\alpha}\rho_1^{\frac{\alpha}{2}}\max\{1,4^{\alpha-1}\rho_1^{-\frac{\alpha-1}{2}}\}.
\end{align}
Then, we can choose $p\in(0,\alpha-1)$ sufficiently small such that
\begin{align}\label{HHp}
\widehat{H}_\alpha(p)\rho_1^{\frac{\alpha}{2}}>L_f+K_2\sqrt{L_g-\rho_1}+K_3(L_g-\rho_1),
\end{align}
 where $\widehat{H}_\alpha(p)$ is given in (\ref{Hal}). Then, for all $y_0\in\mathbb{R}^n$,
\begin{align*}
\mathbb{E}|y(t;y_0)|^p\leq |y_0|^pe^{-\bar{p} t}, \quad \forall t\geq 0.
\end{align*}
where $\bar{p}=p\left(\widehat{H}_\alpha(p)\rho_1^{\frac{\alpha}{2}}-L_f-K_2\sqrt{L_g-\rho_1}-K_3(L_g-\rho_1)\right)$.
\end{theorem}

\begin{proof} If $y_0=0$, the conclusion follows from the trivial solution \(y(t)\equiv0\). In the following, assume that $y_0\ne0$.
Define the stopping times
\begin{align*}
\tau_k=\inf\{t\ge0:|y(t)|\le k^{-1}\},\quad k\ge1.
\end{align*}
Since $|y_0|>0$, for sufficiently large $k$ we have $|y_0|>k^{-1}$. By the right-continuity of $y(t)$ at $0$, it follows that $\tau_k>0$ a.s. for such $k$. Clearly, $\tau_k\le\tau_{k+1}$. Let
\begin{align*}
\tau_\infty:=\lim_{k\to\infty}\tau_k.
\end{align*}
We first prove
\begin{align}\label{tau}
\tau_\infty=\infty,\quad \mathrm{a.s.}
\end{align}
 If it is false, then  $\mathbb P(\tilde{\Omega})>0$, where $\tilde{\Omega}=\{\tau_\infty<\infty\}$. By Theorem \ref{th3.1}, it follows that $y(\tau_\infty)\neq 0$ on $\tilde{\Omega}$. On $\tilde{\Omega}$, if $\tau_k=\tau_\infty$ for infinitely many $k$, then by the definition of $\tau_k$,
\begin{align*}
|y(\tau_\infty)|=|y(\tau_k)|\le k^{-1}
\end{align*}
for infinitely many $k$. Letting $k\rightarrow\infty$ this yields $|y(\tau_\infty)|=0$, which contradicts $y(\tau_\infty)\neq 0$ on $\tilde{\Omega}$. This implies that for all sufficiently large $k$, $\tau_k<\tau_\infty$. We can choose an increasing subsequence $\tau_{k_m}\uparrow\tau_\infty$ as $k_m\rightarrow\infty$. Then,
\begin{align*}
|y(\tau_{k_m})|\le \frac1{k_m}\to 0,
\end{align*}
as $k_m\rightarrow\infty$.
Since $y$ is c\`{a}dl\`{a}g,
\begin{align*}
y(\tau_\infty^-)=\lim_{k_m\to\infty}y(\tau_{k_m})=0.
\end{align*}
Noting that
\begin{align*}
\Delta y(\tau_\infty)=y(\tau_\infty)-y(\tau_\infty-)=g(y(\tau_\infty-))\Delta Z_{\tau_\infty}=0,
\end{align*}
we have $y(\tau_\infty)=y(\tau_\infty-)=0$, which contradicts $y(\tau_\infty)\neq 0$ on $\tilde{\Omega}$. This yields (\ref{tau}) holds. Define the stopping times
\begin{align*}
\sigma_m=\inf\{t\ge0:|y(t)|\ge m\}\wedge m,\qquad m\ge1.
\end{align*}
Since $y(t)$ is c\`{a}dl\`{a}g, for every $T>0$, $\sup_{0\leq t\leq T}|y(t)|<\infty$ a.s. For sufficiently large $m$, both $\inf\{t\geq 0: |y(t)|\geq m\}>T$ and $m>T$ hold, which yields $\sigma_m>T$. Because $T$ is arbitrary,
$\lim\limits_{m\rightarrow\infty}\sigma_m=\infty$, a.s.
 Set
\begin{align*}
\rho_m=\tau_m\wedge\sigma_m.
\end{align*}
 Then
\begin{align}\label{rhom}
\lim\limits_{m\rightarrow\infty}\rho_m=\infty \quad \mathrm{a.s.}
\end{align}

Since $1<\alpha<2$ and $\Gamma(1+s)=s\Gamma(s)$ for $s>0$, we have
\begin{align*}
 \widehat H_\alpha(s) = -\frac{2}{\pi} \Gamma(s)\Gamma(\alpha-s) \sin\frac{\pi s}{2} \cos\frac{\pi(\alpha-s)}{2}.
 \end{align*}
By $\Gamma(s)\sim\frac{1}{s}$, $\sin\frac{\pi s}{2}\sim\frac{\pi s}{2}$, as $s\rightarrow 0$, it follows that
\begin{align*}
\lim\limits_{s\rightarrow 0}\Gamma(s)\sin\frac{\pi s}{2}=\frac{\pi}{2}.
\end{align*}
Consequently,
\begin{align}\label{HG}
\lim_{s\rightarrow0}\widehat H_\alpha(s) = -\Gamma(\alpha)\cos\frac{\alpha\pi}{2}.
\end{align}
Combining (\ref{gamm}) with (\ref{HG}) gives
\begin{align*}
\lim_{s\rightarrow0} \widehat H_\alpha(s)\rho_1^{\alpha/2} >L_f+K_2\sqrt{L_g-\rho_1}+K_3(L_g-\rho_1).
\end{align*}
Set
\begin{align*}
\varepsilon=\frac12\left( -\Gamma(\alpha)\cos\frac{\alpha\pi}{2}\rho_1^{\alpha/2}-L_f-K_2\sqrt{L_g-\rho_1}-K_3(L_g-\rho_1) \right)>0.
\end{align*}
By the definition of the limit, there exists $\delta>0$ such that for every $0<s<\delta$,
\begin{align}\label{liv}
\left| \widehat H_\alpha(s)\rho_1^{\alpha/2} +\Gamma(\alpha)\cos\frac{\alpha\pi}{2}\rho_1^{\alpha/2} \right|<\varepsilon.
\end{align}
By (\ref{gamm}) and (\ref{liv}), we obtain
\begin{align}\label{Halpha}
\widehat H_\alpha(s)\rho_1^{\alpha/2} > -\Gamma(\alpha)\cos\frac{\alpha\pi}{2}\rho_1^{\alpha/2}-\varepsilon>L_f+K_2\sqrt{L_g-\rho_1}+K_3(L_g-\rho_1).
 \end{align}
Let $p_0=\min\left\{\delta,\frac{\alpha-1}{2}\right\}$. Then, for every $p\in(0,p_0)$, we have $p\in(0,\alpha-1)$. Since \eqref{Halpha} holds, taking \(s=p\) yields \eqref{HHp}. Let $V_1(y)=|y|^p$ and $H(t)=e^{\bar{p}t}V_1(y(t))$. Since $V_1\in C^2(\mathbb R^n\setminus\{0\})$ and \(\rho_m\uparrow\infty\) and (\ref{rhom}), the It\^o's formula gives
\begin{align}\label{Htm}
H(t\wedge\rho_m)&= V_1(y_0)-F_m(t)+\widetilde{M}_m(t),
\end{align}
where
\begin{align*}
\mathcal{L}V_1(y)&=p|y|^{p-2} y^\top f(y) + \int_{\mathbb{R}} \Bigl( |y+g(y)z|^{p} - |y|^{p} - p|y|^{p-2} y^\top g(y)z \Bigr) \frac{C_\alpha}{|z|^{\alpha+1}}dz,\nonumber\\
M(t)&=\int_0^{t}\int_{\mathbb R}  \bigl( |y(s^-)+g(y(s^-))z|^{p} - |y(s^-)|^{p} \bigr)\widetilde{N}(ds,dz),\nonumber\\
\widetilde{M}_m(t)&=\int_0^{t\wedge\rho_m} e^{\bar{p}s}dM(s),\nonumber\\
F_m(t)&=-\int_0^{t\wedge\rho_m} e^{\bar{p}s}
\Bigl(\bar{p}V_1(y(s-))+\mathcal{L}V_1(y(s-))\Bigr)ds.
\end{align*}
To study the asymptotic behaviour of the generator term $\mathcal{L}V_1(y)$, fix $y\in\mathbb R^n\setminus\{0\}$, and set $u=y/|y|\in \mathbb{S}^{n-1}$, $\varpi=\frac{g(y)}{|y|}$, $\mu_1=u^\top \varpi$. By Assumption \ref{ass3.2},
\begin{align*}
|\mu_1|^2
=\frac{|y^\top g(y)|^2}{|y|^4}
\geq\rho_1>0.
\end{align*}
It follows from Assumption \ref{ass3.1} and Assumption \ref{ass3.2} that
\begin{align}\label{3.14}
\mathcal{L}V_1(y)&=p|y|^{p-2}y^\top f(y)+|y|^p\int_{\mathbb{R}}(|u+ \varpi z|^{p}-1- p\mu_1z)\frac{C_\alpha}{|z|^{\alpha+1}}dz\nonumber\\
&\leq V_1(y)\left(pL_f+J_1+J_2\right),
\end{align}
where
\begin{align*}
J_1&=\int_{\mathbb{R}}\left(|1+\mu_1z|^{p}-1-p\mu_1z\right)\frac{C_\alpha}{|z|^{\alpha+1}}dz,\\
J_2&=\int_{\mathbb{R}}\left(|u+ \varpi z|^p-|1+\mu_1 z|^p  \right)\frac{C_\alpha}{|z|^{\alpha+1}}dz.
\end{align*}
Define $F(z)=|1+\mu_1z|^p-1-p\mu_1z$. Choose $0<\bar{\delta}< \min\{ 1, \frac{1}{2|\mu_1|}\}$. For $|z|\leq \bar{\delta}$, we have $1+\mu_1z\in[\frac{1}{2},\frac{3}{2}]$. The function $F$ is twice continuously differentiable on $[-\bar{\delta},\bar{\delta}]$, with
\begin{align*}
F(0)=F'(0)=0, \quad F''(z)=p(p-1)\mu_1^2(1+\mu_1z)^{p-2}.
\end{align*}
This implies that
\begin{align*}
|F(z)|\leq p(1-p)\mu_1^22^{1-p}|z|^2,\quad |z|\leq\bar{\delta}.
\end{align*}
We then obtain
\begin{align}\label{F1}
\int_{|z|\leq\bar{\delta}}|F(z)|\frac{C_\alpha}{|z|^{\alpha+1}}dz
&\leq C_\alpha p(1-p)\mu_1^22^{1-p}
\int_{|z|\leq\bar{\delta}}|z|^{1-\alpha}dz\nonumber\\
&=\frac{C_\alpha p(1-p)\mu_1^22^{2-p}}{2-\alpha}
\bar{\delta}^{2-\alpha}.
\end{align}
Since $0<p<\alpha-1<1$,
\begin{align*}
|F(z)|\leq2+|\mu_1|^p|z|^p+p|\mu_1||z|.
\end{align*}
It follows that
\begin{align}\label{F2}
\int_{|z|>\bar{\delta}}|F(z)|\frac{C_\alpha}{|z|^{\alpha+1}}dz
&\leq2C_\alpha\int_{\bar{\delta}}^\infty
\left(2+|\mu_1|^pr^p+p|\mu_1|r\right)r^{-1-\alpha}dr\nonumber\\
&=2C_\alpha\left(
\frac{2\bar{\delta}^{-\alpha}}{\alpha}
+\frac{|\mu_1|^p\bar{\delta}^{p-\alpha}}{\alpha-p}
+\frac{p|\mu_1|\bar{\delta}^{1-\alpha}}{\alpha-1}
\right).
\end{align}
Combining (\ref{F1}) with (\ref{F2}) gives that
\begin{align}\label{F3}
\int_{\mathbb R}|F(z)|\frac{C_\alpha}{|z|^{\alpha+1}}dz\leq\frac{C_\alpha p(1-p)\mu_1^22^{2-p}}{2-\alpha}
\bar{\delta}^{2-\alpha}+2C_\alpha\left(
\frac{2\bar{\delta}^{-\alpha}}{\alpha}
+\frac{|\mu_1|^p\bar{\delta}^{p-\alpha}}{\alpha-p}
+\frac{p|\mu_1|\bar{\delta}^{1-\alpha}}{\alpha-1}
\right)<\infty.
\end{align}
For each $t>0$, let $q_t$ denote the continuous probability density function of the random variable $Z_t$. By the self-similarity of the symmetric $\alpha$-stable process,
\begin{align}\label{qt}
q_t(z)=t^{-\frac{1}{\alpha}}q_1(t^{-\frac{1}{\alpha}}z),\quad t>0,\,\, z\in \mathbb{R},
\end{align}
where $q_1$ is the probability density function of $Z_1$. By Proposition 7.1.2 in Kolokol'tsov \cite{Kolokoltsov2011} and the symmetry of $q_1$, we obtain
\begin{align}\label{q1}
q_1(x)\sim \frac{\Gamma(1+\alpha)\sin \left(\frac{\pi\alpha}{2}\right) }{\pi} |x|^{-1-\alpha}, \quad |x|\rightarrow \infty.
\end{align}
Using the reflection and duplication formulas
\begin{align*}
\Gamma\left(\frac{\alpha}{2}\right)\Gamma\left(1-\frac{\alpha}{2}\right)
=\frac{\pi}{\sin\left(\frac{\pi\alpha}{2}\right)},
\qquad
\Gamma\left(\frac{\alpha}{2}\right)\Gamma\left(\frac{1+\alpha}{2}\right)
=2^{1-\alpha}\sqrt{\pi}\Gamma(\alpha),
\end{align*}
together with \(\Gamma(1+\alpha)=\alpha\Gamma(\alpha)\), we obtain
\begin{align}\label{Calpha}
\frac{\Gamma(1+\alpha)\sin\left(\frac{\pi\alpha}{2}\right)}{\pi}=C_\alpha.
\end{align}
Define $h(x)=|x|^{1+\alpha}q_1(x)$ for $x\in\mathbb R$. Since $q_1$ is continuous, $h$ is continuous on $\mathbb R$. Moreover, by (\ref{q1}) and (\ref{Calpha}), it follows that
\begin{align}\label{hx}
\lim_{|x|\to\infty}h(x)=C_\alpha.
 \end{align}
It follows that $h$ is bounded, and there exists a constant $K_\alpha>0$ such that
\begin{align}\label{Kal}
0\leq h(x)\leq K_\alpha,\qquad x\in\mathbb R.
\end{align}
By (\ref{qt}) and (\ref{Kal}), we have
\begin{align}\label{Ka2}
\frac{q_t(z)}{t}=|z|^{-1-\alpha}h(t^{-1/\alpha}z)\leq K_\alpha|z|^{-1-\alpha},\quad t>0,\, z\neq0.
\end{align}
The function $K_\alpha|F(z)||z|^{-1-\alpha}$ is integrable on $\mathbb R$ by (\ref{F3}), and
\begin{align*}
\left|F(z)\frac{q_t(z)}{t}\right|
\leq K_\alpha|F(z)||z|^{-1-\alpha},
\qquad t>0,\ z\neq0.
\end{align*}
Moreover, for each fixed $z\neq0$, (\ref{qt}) and (\ref{hx}) give
\begin{align}\label{Ka4}
\lim\limits_{t\rightarrow 0}\frac{q_t(z)}{t}=\lim\limits_{t\rightarrow 0}|z|^{-1-\alpha}h(t^{-1/\alpha}z)=|z|^{-1-\alpha}C_\alpha,\ z\neq0.
\end{align}
Since $\mathbb E[Z_t]=0$, the dominated convergence theorem and (\ref{Ka4}) yield
\begin{align}\label{3.15}
\lim_{t\rightarrow0}\frac{\mathbb{E}|1+\mu_1Z_t|^p-1}{t}
&=\lim_{t\rightarrow0}\frac{\mathbb{E}[F(Z_t)]}{t}\nonumber\\
&=\lim_{t\rightarrow0}\int_{\mathbb{R}}F(z)\frac{q_t(z)}tdz\nonumber\\
&=\int_{\mathbb{R}} F(z)\frac{C_\alpha}{|z|^{1+\alpha}}dz\nonumber\\
&=J_1.
\end{align}
By Lemma \ref{le2.1},
 \begin{align}\label{3.16}
 \mathbb{E}|1+\mu_1Z_t|^p=\frac{2\Gamma(1+p)\sin\frac{\pi p}{2}}{\pi}\int_0^\infty\frac{1-e^{-|\mu_1|^\alpha tr^\alpha}\cos r}{r^{p+1}}dr.
\end{align}
Applying Lemma \ref{le2.4} with $\vartheta=\mu_1$, and using (\ref{3.15}) and (\ref{3.16}) we obtain
\begin{align*}
J_1=-p\widehat{H}_\alpha(p)|\mu_1|^\alpha,
\end{align*}
where $\widehat{H}_\alpha(p)$ is given in (\ref{Hal}). Since $\widehat H_\alpha(p)>0$, Assumption \ref{ass3.2} gives
\begin{align}\label{3.19}
J_1\leq -p\widehat{H}_\alpha(p)\rho_1^{\frac{\alpha}{2}}.
\end{align}

 Let $\widetilde{\Omega}_1=\{|z|\leq \delta_1\}$, $\widetilde{\Omega}_2=\{\delta_1<|z|\leq R_1, |z+\mu_1^{-1}|\geq \varepsilon_1\}$, $\widetilde{\Omega}_3=\{\delta_1<|z|\leq R_1, |z+\mu_1^{-1}|< \varepsilon_1\}$, and $\widetilde{\Omega}_4=\{ |z|>R_1\}$,
where $\delta_1=\frac{1}{4\sqrt{L_g}}$, $R_1=\frac{4}{\sqrt{\rho_1}}$, and $\varepsilon_1=\frac{1}{4}\min\{1/\sqrt{L_g},1/\sqrt{\rho_1}\}$. Using $|u+z\varpi|^2=(1+\mu_1z)^2+(\mu_2-\mu_1^2)z^2$,  we decompose $J_2$ as
\begin{align}\label{3.20}
J_2=\sum_{i=1}^{4}J_2^{(i)},
\end{align}
where
\begin{align*}
J_2^{(i)}=C_\alpha\int_{\widetilde{\Omega}_i} \left[((1+\mu_1z)^2+(\mu_2-\mu_1^2)z^2)^{p/2}-|1+\mu_1z|^p\right] \frac{1}{|z|^{\alpha+1}}dz.
\end{align*}
By
\begin{align}\label{a1a2}
0\leq (a_1^2+a_2^2)^{p/2}-a_1^p\leq \frac{p}{2}a_2^2a_1^{p-2}, \quad a_1>0, a_2\geq 0,
\end{align}
it follows that
\begin{align}\label{3.21}
J_2^{(1)}&\leq C_\alpha \frac{p}{2}(\mu_2-\mu_1^2)\int_{\widetilde{\Omega}_1} \frac{1}{|1+\mu_1z|^{2-p}}|z|^{1-\alpha}dz\nonumber\\
&\leq p(\mu_2-\mu_1^2)\frac{8}{9}C_\alpha \int_{|z|\leq \delta_1}|z|^{1-\alpha}dz\nonumber\\
&\leq p (\mu_2-\mu_1^2) \frac{16C_\alpha}{9(2-\alpha)}4^{\alpha-2}L_g^{\frac{\alpha-2}{2}},
\end{align}
\begin{align}\label{3.22}
J_2^{(2)}&\leq C_\alpha \int_{\widetilde{\Omega}_2}\frac{p}{2}(\mu_2-\mu_1^2)z^2|1+\mu_1z|^{p-2}\frac{1}{|z|^{\alpha+1}}dz\nonumber\\
&\leq C_\alpha \int_{\widetilde{\Omega}_2}\frac{p}{2}(\mu_2-\mu_1^2)R_1^2(\sqrt{\rho_1}\varepsilon_1)^{p-2}\frac{1}{|z|^{\alpha+1}}dz\nonumber\\
&\leq C_\alpha\frac{p}{2}(\mu_2-\mu_1^2)R_1^2(\sqrt{\rho_1}\varepsilon_1)^{-2}\int_{\delta_1<|z|\leq R_1}|z|^{-1-\alpha}dz\nonumber\\
&\leq p(\mu_2-\mu_1^2)C_\alpha \frac{256}{\alpha\rho_1^2}\max\{L_g,\rho_1\}\left(4^\alpha L_g^{\frac{\alpha}{2}}-4^{-\alpha}\rho_1^{\frac{\alpha}{2}}\right),
\end{align}
and
\begin{align}\label{3.23}
J_2^{(4)}&\leq C_\alpha \int_{\widetilde{\Omega}_4}\frac{p}{2}(\mu_2-\mu_1^2)z^2|1+\mu_1z|^{p-2}\frac{1}{|z|^{\alpha+1}}dz\nonumber\\
&\leq  C_\alpha \int_{\widetilde{\Omega}_4}\frac{p}{2}(\mu_2-\mu_1^2)z^2\left|\frac{1}{2}\sqrt{\rho_1}|z|\right|^{p-2}\frac{1}{|z|^{\alpha+1}}dz\nonumber\\
&\leq pC_\alpha(\mu_2-\mu_1^2)\left(\frac{2}{\sqrt{\rho_1}}\right)^{2-p}\frac{R_1^{p-\alpha}}{\alpha-p}\nonumber\\
&\leq pC_\alpha (\mu_2-\mu_1^2)\max\left\{\frac{4}{\rho_1},\frac{2^{3-\alpha}}{\rho_1^{\frac{3-\alpha}{2}}}\right\} 4^{-\alpha}\rho_1^{\frac{\alpha}{2}}\max\left\{1,4^{\alpha-1}\rho_1^{-\frac{\alpha-1}{2}}\right\}.
\end{align}
If $\mu_2-\mu_1^2=0$, then $J_2^{(3)}(u)=0$. It remains to consider $\mu_2-\mu_1^2>0$. By (\ref{a1a2}), it follows that
\begin{align}\label{J23}
J_2^{(3)}&\leq \left(\frac{4}{3}\right)^{1+\alpha}C_\alpha |\mu_1|^{1+\alpha}\int_{|s|<\varepsilon_1} \left[(\mu_1^2s^2+(\mu_2-\mu_1^2)(s-\mu_1^{-1})^2)^{p/2}-|\mu_1|^p|s|^p  \right] ds\nonumber\\
&\leq \left(\frac{4}{3}\right)^{1+\alpha}C_\alpha |\mu_1|^{1+\alpha}(J_2^A+J_2^B),
\end{align}
where
\begin{align*}
J_2^A&=\int_{\{|s|<\varepsilon_1,|\mu_1s|\geq \sqrt{\mu_2-\mu_1^2}|s-\mu_1^{-1}|\}} \left[(\mu_1^2s^2+(\mu_2-\mu_1^2)(s-\mu_1^{-1})^2)^{p/2}-|\mu_1|^p|s|^p  \right] ds,\nonumber\\
J_2^B&=\int_{\{|s|<\varepsilon_1,|\mu_1s|< \sqrt{\mu_2-\mu_1^2}|s-\mu_1^{-1}|\}} \left[(\mu_1^2s^2+(\mu_2-\mu_1^2)(s-\mu_1^{-1})^2)^{p/2}-|\mu_1|^p|s|^p  \right] ds.
\end{align*}
By (\ref{a1a2}), we have
\begin{align}\label{J2A}
J_2^A&\leq p \frac{25}{32}(\mu_2-\mu_1^2)|\mu_1|^{p-4}\int_{\frac{3}{4}\frac{\sqrt{\mu_2-\mu_1^2}}{\mu_1^2}<|s|<\varepsilon_1}|s|^{p-2}ds\nonumber\\
&\leq p\frac{25}{16}\frac{3^{p-1}}{4^{p-1}}\frac{1}{1-p}|\mu_1|^{-p-2}(\mu_2-\mu_1^2)^{\frac{p+1}{2}}.
\end{align}
Using the inequality
\begin{align*}
\frac{(1+r^2)^{p/2}-r^p}{p}\leq 2(1+|\log r|),\quad 0<r<1,
\end{align*}
it follows that
\begin{align}\label{J2B}
J_2^B
&=\tilde{\gamma}
\int_{D_B}
\left[
\left(\mu_1^2\tilde{\gamma}^2s^2
+(\mu_2-\mu_1^2)(\tilde{\gamma}s-\mu_1^{-1})^2\right)^{p/2}
-|\mu_1|^p|\tilde{\gamma}s|^p
\right]ds
\nonumber\\
&\leq
4p\max\left\{
1,(\sqrt{L_g}R_1)^{\frac{\alpha-1}{2}}
\right\}\tilde{\gamma}
\int_{D_B}
(1+|\log|s||)\,ds
\nonumber\\
&\leq
4p\max\left\{
1,(\sqrt{L_g}R_1)^{\frac{\alpha-1}{2}}
\right\}\tilde{\gamma}
\int_{\left\{|s|\leq\frac54\right\}}
(1+|\log|s||)\,ds
\nonumber\\
&=
4p\max\left\{
1,(\sqrt{L_g}R_1)^{\frac{\alpha-1}{2}}
\right\}\tilde{\gamma}
\left(4+\frac52\log\frac54\right),
\end{align}
where $\tilde{\gamma}=\frac{\sqrt{\mu_2-\mu_1^2}}{\mu_1^2}$ and $D_B=\{s\in\mathbb{R}: |\tilde{\gamma}s|<\varepsilon_1, |s|<|\mu_1||\tilde{\gamma}s-\mu_1^{-1}| \}$. By (\ref{J23})-(\ref{J2B}), we get that
\begin{align}\label{3.24}
J_2^{(3)}\leq p\sqrt{\mu_2-\mu_1^2}K_2,
\end{align}
where $K_2$ is given in (\ref{K2}). Applying $\mu_2-\mu_1^2\leq L_g-\rho_1$, (\ref{3.20}), (\ref{3.21})-(\ref{3.23}), and (\ref{3.24}), we can conclude that
\begin{align}\label{3.27}
J_2&\leq p[K_2\sqrt{L_g-\rho_1}+K_3(L_g-\rho_1)],
\end{align}
where $K_3$ is given in (\ref{K3}). By (\ref{HHp}), (\ref{3.14}), (\ref{3.19}), and (\ref{3.27}), we have
\begin{align}\label{3.28}
\mathcal{L}V_1(y)\leq -V_1(y)\bar{p}.
\end{align}
Choose $\delta_m=\min\left\{1,\frac1{2m^2\sqrt{L_g}}\right\}$. Note that
\begin{align}\label{Theta}
\mathbb{E}\left[(M(t\wedge\rho_m))^2\right] &=C_\alpha\mathbb{E}\left[\int_0^{t\wedge\rho_m}\int_{\mathbb R}  \bigl(|y(s^-)+g(y(s^-))z|^{p} - |y(s^-)|^{p} \bigr)^2\frac{1}{|z|^{\alpha+1}}dzds\right]\nonumber\\
&=:\Theta_1+\Theta_2+\Theta_3,
\end{align}
where
\begin{align*}
\Theta_1&=C_\alpha\mathbb{E}\left[\int_0^{t\wedge\rho_m}\int_{|z|\leq \delta_m}  \bigl(|y(s^-)+g(y(s^-))z|^{p} - |y(s^-)|^{p} \bigr)^2\frac{1}{|z|^{\alpha+1}}dzds\right],\nonumber\\
\Theta_2&=C_\alpha\mathbb{E}\left[\int_0^{t\wedge\rho_m}\int_{\delta_m<|z|\leq 1}  \bigl(|y(s^-)+g(y(s^-))z|^{p} - |y(s^-)|^{p} \bigr)^2\frac{1}{|z|^{\alpha+1}}dzds\right],\nonumber\\
\Theta_3&=C_\alpha\mathbb{E}\left[\int_0^{t\wedge\rho_m}\int_{|z|>1}  \bigl(|y(s^-)+g(y(s^-))z|^{p} - |y(s^-)|^{p} \bigr)^2\frac{1}{|z|^{\alpha+1}}dzds\right].
\end{align*}
Note that
\begin{align}\label{ym}
m^{-1}\leq |y(s^-)|\leq m, \quad 0<s\leq \rho_m.
\end{align}
For $\Theta_1$, applying Assumption \ref{ass3.2}, (\ref{ym}) and the mean value theorem, we get that there exists $\theta\in(0,1)$ such that
\begin{align}\label{Theta1}
\Theta_1&\leq C_\alpha\mathbb{E}\left[\int_0^{t\wedge\rho_m}\int_{|z|\leq \delta_m}  \left(p|y(s^-)+\theta g(y(s^-))z|^{p-1}|g(y(s^-))z|\right)^2\frac{1}{|z|^{\alpha+1}}dzds\right]\nonumber\\
&\leq  C_\alpha \mathbb{E}\left[\int_0^{t\wedge\rho_m}\int_{|z|\leq \delta_m}  \left(p(2m)^{1-p}\sqrt{L_g}m|z|\right)^2\frac{1}{|z|^{\alpha+1}}dzds\right]\nonumber\\
&\leq C_\alpha p^2L_g2^{2-2p}m^{4-2p}\int_0^{t}\int_{|z|\leq 1}  \frac{1}{|z|^{\alpha-1}}dzds\nonumber\\
&\leq C_\alpha p^2L_g2^{2-2p}m^{4-2p}\frac{2}{2-\alpha}t.
\end{align}
For $\Theta_2$, by (\ref{ym}) and Assumption \ref{ass3.2}, it follows that
\begin{align}\label{Theta2}
\Theta_2&\leq C_\alpha\mathbb{E}\left[ \int_0^{t\wedge\rho_m}\int_{\delta_m<|z|\leq 1}(\sqrt{L_g}m)^{2p}|z|^{2p-\alpha-1}dzds  \right]\nonumber\\
&=C_\alpha(\sqrt{L_g}m)^{2p}\frac{2}{\alpha-2p}(\delta_m^{2p-\alpha}-1)\mathbb{E}[t\wedge\rho_m]\nonumber\\
&\leq C_\alpha(\sqrt{L_g}m)^{2p}\frac{2}{\alpha-2p}(\delta_m^{2p-\alpha}-1)t.
\end{align}
For $\Theta_3$, combining (\ref{ym}) with Assumption \ref{ass3.2} gives that
\begin{align}\label{Theta3}
\Theta_3&\leq C_\alpha\mathbb{E}\left[\int_0^{t\wedge\rho_m}\int_{|z|> 1}(\sqrt{L_g}m)^{2p}|z|^{2p-\alpha-1}dzds  \right]\nonumber\\
&=C_\alpha (\sqrt{L_g}m)^{2p}\frac{2}{\alpha-2p}\mathbb{E}[t\wedge\rho_m]\nonumber\\
&\leq C_\alpha (\sqrt{L_g}m)^{2p}\frac{2}{\alpha-2p}t.
\end{align}
By (\ref{Theta}), and (\ref{Theta1})-(\ref{Theta3}), we have
\begin{align*}
\mathbb{E}\left[(M(t\wedge\rho_m))^2\right]\leq C_\alpha p^2L_g2^{2-2p}m^{4-2p}\frac{2}{2-\alpha}t+C_\alpha(\sqrt{L_g}m)^{2p}\frac{2}{\alpha-2p}(\delta_m^{2p-\alpha})t<\infty.
\end{align*}
Consequently, $M(t\wedge\rho_m)$ is a square-integrable martingale. By \eqref{rhom}, we get $\{M(t)\}_{t\geq 0}$ is a local martingale. Hence, $\{\widetilde{M}_m(t)\}_{t\geq 0}$ is a local martingale. Using (\ref{3.28}),  we obtain that $F_m(t)$ is a nonnegative nondecreasing process. It then follows from \eqref{Htm} that $\{H(t\wedge\rho_m)\}_{t\ge0}$ is a nonnegative local supermartingale. Hence, $\{H(t\wedge\rho_m)\}_{t\ge0}$ is a true supermartingale and
\begin{align*}
\mathbb{E}[H(t\wedge\rho_m)\mid\mathcal F_s]\le H(s\wedge\rho_m),\qquad 0\le s\le t.
\end{align*}
Applying the conditional form of Fatou's lemma and \eqref{rhom}, we have
\begin{align*}
\mathbb{E}[H(t)\mid\mathcal F_s]
\le \liminf_{m\to\infty}\mathbb E[H(t\wedge\rho_m)\mid\mathcal F_s]
\le \liminf_{m\to\infty}H(s\wedge\rho_m)
=H(s).
\end{align*}
It follows that $\{H(t)\}_{t\geq 0}$ is a right-continuous nonnegative supermartingale, and
\begin{align}\label{H}
\mathbb{E}[H(t)]\le \mathbb E[H(0)]=|y_0|^p.
\end{align}
Consequently,
\begin{align*}
\mathbb{E}|y(t)|^p\leq |y_0|^pe^{-\bar{p}t},\quad t\geq 0.
\end{align*}
\end{proof}

\section{Almost sure exponential stability of solutions to MSS (\ref{model})}
In this section, we establish the almost sure exponential stability for the MSS (\ref{model}).

\begin{theorem}\label{thm5.1} Suppose that the assumptions of Theorem \ref{th4.1} are satisfied. Let \(p\in(0,\alpha-1)\) and \(\bar p>0\) be chosen as in Theorem \ref{th4.1}. Then for any initial value $y_0\neq0$, the solution $\{y(t)\}_{t\ge0}$ of equation (\ref{model}) satisfies
\begin{equation}\label{expstable}
\limsup_{t\to\infty}\frac{1}{t}\log|y(t)| \le -\left(\widehat{H}_\alpha(p)\rho_1^{\frac{\alpha}{2}}-L_f-K_2\sqrt{L_g-\rho_1}-K_3(L_g-\rho_1)\right) \qquad \text{a.s.}
\end{equation}
\end{theorem}

\begin{proof}
For any fixed $\lambda\in(0,\bar{p}/p)$ and each integer $n\ge0$, define
\begin{align*}
A_n=\Bigl\{\sup_{t\in[n,n+1]}|y(t)|^p\ge e^{-p\lambda n}\Bigr\},
\end{align*}
 On $A_n$,
\begin{align*}
\sup_{t\in[n,n+1]}H(t)
\ge e^{\bar{p}n}\sup_{t\in[n,n+1]}|y(t)|^p
\ge e^{(\bar{p}-p\lambda)n}.
\end{align*}
As a result,
\begin{align*}
A_n\subseteq \left\{\sup_{t\in[n,n+1]}H(t)\ge e^{(\bar{p}-p\lambda)n}\right\}.
\end{align*}
 Doob's maximal inequality and (\ref{H}) give that
\begin{align*}
\mathbb{P}(A_n) &\le \mathbb{P}\left(\sup_{t\in[n,n+1]}H(t)\ge e^{(\bar{p}-p\lambda)n}\right)\nonumber\\
&\leq \frac{\mathbb{E}[H(n)]}{e^{(\bar{p}-p\lambda)n}}\nonumber\\
&\leq |y_0|^p e^{-(\bar{p}-p\lambda)n}.
\end{align*}
By the Borel-Cantelli lemma, there exists an almost surely finite integer \(N(\omega)\) such that for all $n\ge N(\omega)$,
\begin{align*}
\sup_{t\in[n,n+1]}|y(t)|\le e^{-\lambda n}.
\end{align*}
For any $t\ge N(\omega)+1$, taking $n=\lfloor t\rfloor$, we get
\begin{align*}
|y(t)|\le e^{-\lambda n}\le e^{\lambda}e^{-\lambda t}.
\end{align*}
Thus,
\begin{align*}
\limsup_{t\to\infty}\frac1t\log|y(t)|\le -\lambda.
\end{align*}
Since $\lambda\in(0,\bar{p}/p)$ is arbitrary, letting $\lambda\rightarrow \bar{p}/p$ along a countable sequence yields \eqref{expstable} holds. This completes the proof.
\end{proof}

\begin{remark}
Let $n=1$ and consider the scalar SDE
\begin{align}
dy(t)=\theta y(t)dt+b y(t)dZ_t,\quad y_0\in\mathbb{R},
\end{align}
where $\theta\in\mathbb{R}$, $b\in\mathbb{R}\setminus\{0\}$. Then Assumptions \ref{ass3.1} and \ref{ass3.2} are satisfied with
\begin{align}
L_f=|\theta|,\quad L_g=\rho_1=b^2.
\end{align}
Hence, the general condition \eqref{gamm} reduces to
\begin{align}\label{1dcond}
-\Gamma(\alpha)\cos\frac{\alpha\pi}{2}|b|^\alpha>|\theta|.
\end{align}
If (\ref{1dcond}) holds, we can choose $p\in(0,\alpha-1)$ sufficiently small such that
$p\left(\widehat H_\alpha(p)|b|^\alpha-|\theta|\right)>0$, where $\widehat{H}_\alpha(p)$ is given in (\ref{Hal}). Then, for all $y_0\in\mathbb{R}$,
\begin{align*}
\mathbb{E}|y(t; y_0)|^p\le|y_0|^p e^{-p\left(\widehat H_\alpha(p)|b|^\alpha-|\theta|\right) t},\qquad t\ge0.
\end{align*}
Moreover, for any initial value $y_0\neq0$,
\begin{align*}
\limsup_{t\to\infty}\frac1t\log|y(t)|
\le|\theta|-\widehat H_\alpha(p)|b|^\alpha<0
\qquad\text{a.s.}
\end{align*}
\end{remark}

\section{Numerical example}
In this section, we present a numerical example illustrating the theoretical result. Consider the two-dimensional stochastic system
\begin{align}\label{ex}
dy(t) = By(t)dt + Ay(t)dZ(t),
\end{align}
where
\begin{align*}
A = \begin{pmatrix} 7.7 & 0.1 \\ -0.1 & 7.7 \end{pmatrix}, \qquad B = \begin{pmatrix} 0.5 & 0.1 \\ 0.05 & 0.45 \end{pmatrix},
\qquad
y(0) = \begin{pmatrix} 1 \\ -0.5 \end{pmatrix},
\end{align*}
and $\{Z(t)\}_{t\ge 0}$ is a standard symmetric $\alpha$-stable process with stability index $\alpha = 1.5$. For the above matrix $A$, we have
\begin{align*}
|Ay|^2 \le L_g |y|^2,\qquad |y^\top A y|^2 \ge \rho_1 |y|^4.
\end{align*}
We may take
\begin{align*}
L_g= ||A||^2 = 59.3,\qquad
\rho_1 = \left(\lambda_{\min}\left(\tfrac{A+A^\top}{2}\right)\right)^2 = 59.29.
\end{align*}
Then, the Lipschitz constant can be taken as
\begin{align*}
L_f = ||B|| = 0.5547.
\end{align*}
With $\alpha = 1.5$, We have $C_\alpha \approx 0.2992$, $K_2\approx 81.8590387084$, $K_3 \approx 145.2595$. This implies that
\begin{align*}
-\Gamma(\alpha)\cos\frac{\alpha\pi}{2}\rho_1^{\alpha/2}
-L_f-K_2\sqrt{L_g-\rho_1}-K_3(L_g-\rho_1)
\approx 3.1963>0,
\end{align*}
i.e., condition \eqref{gamm} holds. Next, taking $p=0.01$, we compute $\widehat H_\alpha(p)\approx 0.6130$, so that
\begin{align*}
\widehat H_\alpha(p)\rho_1^{\alpha/2}
-L_f-K_2\sqrt{L_g-\rho_1}-K_3(L_g-\rho_1)
\approx 2.9049>0.
\end{align*}
Therefore, by Theorem \ref{th4.1}, the noise-perturbed system \eqref{ex} is exponentially stable in the $p$-th moment. Moreover, by Theorem \ref{thm5.1}, it is almost surely exponentially stable.

Let $T>0$ be the observation time horizon and $n$ be the number of time steps, with $\Delta=T/n$. We approximate \eqref{ex} via the Euler scheme (see, e.g., Jacod \cite{Jacod2004} and Jacod et al. \cite{Jacod2005}):
\begin{align*}
y_{t_{k+1}}=y_{t_k}+B\,y_{t_k}\,\Delta+A\,y_{t_k}\,\Delta Z_{t_k},
\end{align*}
where $t_k=\Delta k$ and $\Delta Z_{t_k}=Z_{t_{k+1}}-Z_{t_k}\overset{d}{=}\Delta^{1/\alpha}\widetilde{Z}_k$, $k=0,1,\ldots,n-1$. The i.i.d. random variables $\widetilde{Z}_k$ are generated according to the algorithm of Janicki and Weron \cite{Janicki}. Let $\{(V_k,W_k)\}_{k\ge0}$ be independent and identically distributed pairs, where $V_k$ is uniform on $(-\pi/2,\pi/2)$, $W_k$ is exponentially distributed with mean one, and $V_k$ and $W_k$ are independent. Then the corresponding symmetric $\alpha$-stable random variable is given by
\begin{align*}
\widetilde{Z}_k=\frac{\sin(\alpha V_k)}{(\cos V_k)^{1/\alpha}}\left(\frac{\cos((1-\alpha)V_k)}{W_k}\right)^{(1-\alpha)/\alpha}.
\end{align*}

 Using the above Euler scheme, we simulate the two-dimensional system \eqref{ex}. Figures~1(a) and 1(b) show representative sample paths of $y_1(t)$ and $y_2(t)$, respectively. Figures~2(a) and 2(b) display the corresponding sample paths of $\ln|y_1(t)|/t$ and $\ln|y_2(t)|/t$, respectively.
\begin{figure}[htbp]
    \centering
    \includegraphics[width=0.7\textwidth]{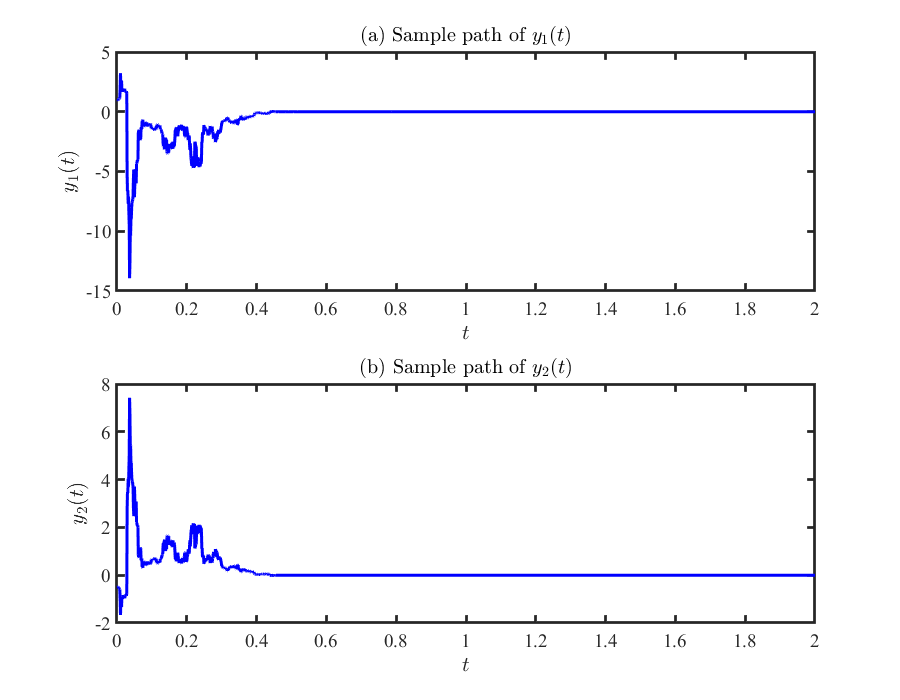}
    \caption{Sample paths of $y_1(t)$ and $y_2(t)$ for system \eqref{ex}.}
    \label{fig:states}
\end{figure}
\begin{figure}[htbp]
    \centering
    \includegraphics[width=0.7\textwidth]{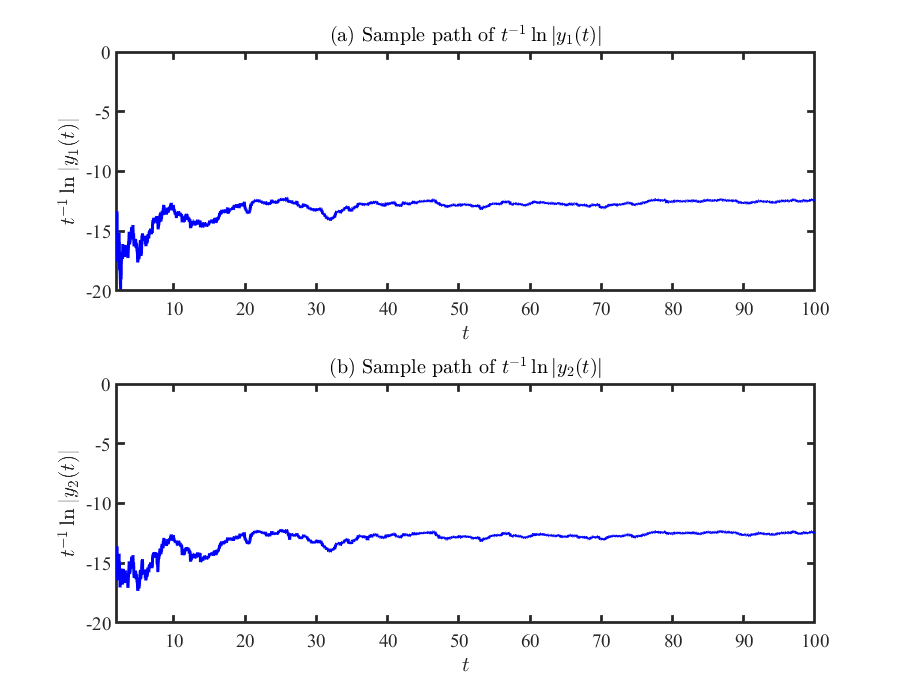}
    \caption{Sample paths of $t^{-1}\ln|y_1(t)|$ and $t^{-1}\ln|y_2(t)|$ for system \eqref{ex}.}
    \label{fig:logrates}
\end{figure}

\section*{Acknowledgments}

The authors are very grateful to the referees and the editor for their insightful and valuable comments which have greatly improved the presentation of the article. This work was supported by the Anhui Provincial Natural Science Foundation (2508085MA016).

\section*{Data availability}

All data generated or analysed during this study are included in this article.

\section*{Disclosure statement}
No potential conflict of interest was reported by the authors.



\section*{Disclosure statement}
No potential conflict of interest was reported by the authors.

\end{document}